\documentclass[12pt]{amsart}

\usepackage{amsmath}
\usepackage{amssymb,mathrsfs}
\usepackage[left=3.5cm,right=3.5cm,top=3.5cm,bottom=3.5cm]{geometry}
\usepackage{graphicx}

\begin{document}

\title{On existence of a triangle with prescribed bisector lengths}

\author{S. F. Osinkin }
\address{ Moscow State University, Russia}
\email[S. F. Osinkin]{tatyana141100@yandex.ru}


\begin{abstract}
We suggest a geometric visualization of the process of constructing a triangle with prescribed
bisectors that makes the existence of such a triangle geometrically evident.
\end{abstract}


\maketitle

{\it Introduction.} --- Since antiquity, in construction problems of elementary geometry
main attention has been paid to possibility
or impossibility of such constructions with only ruler and compass, and problems
of constructing regular $n$-gons,
trisection of angle or squaring of a circle has played an important role in the
development of mathematics. Existence
and number of solutions of these construction problems when one is not confined
to using ruler and compass only is often evident
and as a rule is not discussed in much detail. However, situation changes when
one turns to other problems such as, for example,
problems of constructing triangles with prescribed its three elements.
In particular, a triangle with given lengths of its
sides can be constructed if the `triangle inequalities' are fulfilled and then
the triangle is unique up to an isometry.
In a similar way, the problems of constructing a triangle with given lengths
of its medians or altitudes can be
considered and these constructions with ruler and compass are possible under
conditions that certain inequalities
are fulfilled. Interestingly enough, the problem of constructing a triangle
with prescribed lengths of its bisectors
is very different: it was shown that such a construction with only ruler and
compass is impossible (see, e.g.,
\cite{manin}). Then the question arises: Given three lengths $l_a$, $l_b$, $l_c$,
does there exist a triangle with
its angle bisectors equal to these lengths? Positive answer to this question
was given by Mironescu and Panaitopol
in \cite{mp-1994} by analytical method with the use of the Brower fixed point
theorem. More elementary analytical
proof was suggested by Zhukov and Akulich in \cite{za-2003}. Being absolutely
strict, these proofs lack, in our opinion,
geometric clearness. In this note we suggest such a geometric visualization
of the process of
constructing a triangle with prescribed bisectors that makes the possibility of such a construction practically
evident. Thus, our aim in this note is to present an elementary visual proof of the
Mironescu and Panaitopol theorem that reads:

\smallskip

{\bf Theorem. For any given segments with lengths $l_a$, $l_b$, $l_c$ there exists a triangle
whose bisectors are equal to these prescribed lengths.}

\smallskip

{\it Informal discussion.} --- Before going to the strict proof of this theorem,
we give an informal description
of the idea of our proof in which
the process of constructing a triangle with prescribed bisector lengths is
reduced to two transformations.
Let, for definiteness, the prescribed lengths of bisectors satisfy
inequalities $l_c<l_a<l_b$. We start from
an equilateral triangle with three equal bisector lengths equal to $l_b$.
The first transformation is the decreasing
of the angle $B$ in such a way that the triangle becomes isosceles with
decreasing equal sides $a=c$ and constant
bisector $l_b$. It is evident that bisectors $l_a$ and $l_c$ decrease
also remaining all the time equal to each other.
As a result of this first transformation we can reach such an isosceles
triangle that its two changing bisectors
become equal to the prescribed value $l_a$. Further decrease of the angle
$B$ would lead to decrease of $l_a$,
therefore we combine change of $B$ with rotation of the side $b$ around
the intersection point of $l_b$ and $b$
in such a way that $l_a$ remains fixed (as well as $l_b$) and $l_c$
decreases reaching finally the prescribed
value. Possibility of such transformation follows from geometrically evident observation
that if angles $A$ and $B$ decrease under condition that $l_a$ and $l_b$ remain constant,
then the third angle $C$ increases whereas its bisector $l_c$ decreases. In the limit
$\angle A\to0$, $\angle B\to 0$ we have $l_c\to0$ and hence $l_c$ can reach any prescribed
value $l_c<l_a,\,l_b$.
This visualization of the construction makes intuitively plausible that a triangle with prescribed
bisector lengths $l_c<l_a<l_b$ does exist.

\smallskip

{\it Preliminary lemmas.} --- Now we turn to formal realization of the above idea.
As a preliminary step, we shall prove several Lemmas about angles and sides of a triangle,
and its bisectors.

\begin{figure}[h]
\begin{center}
\includegraphics[width=6cm,height=6cm,clip]{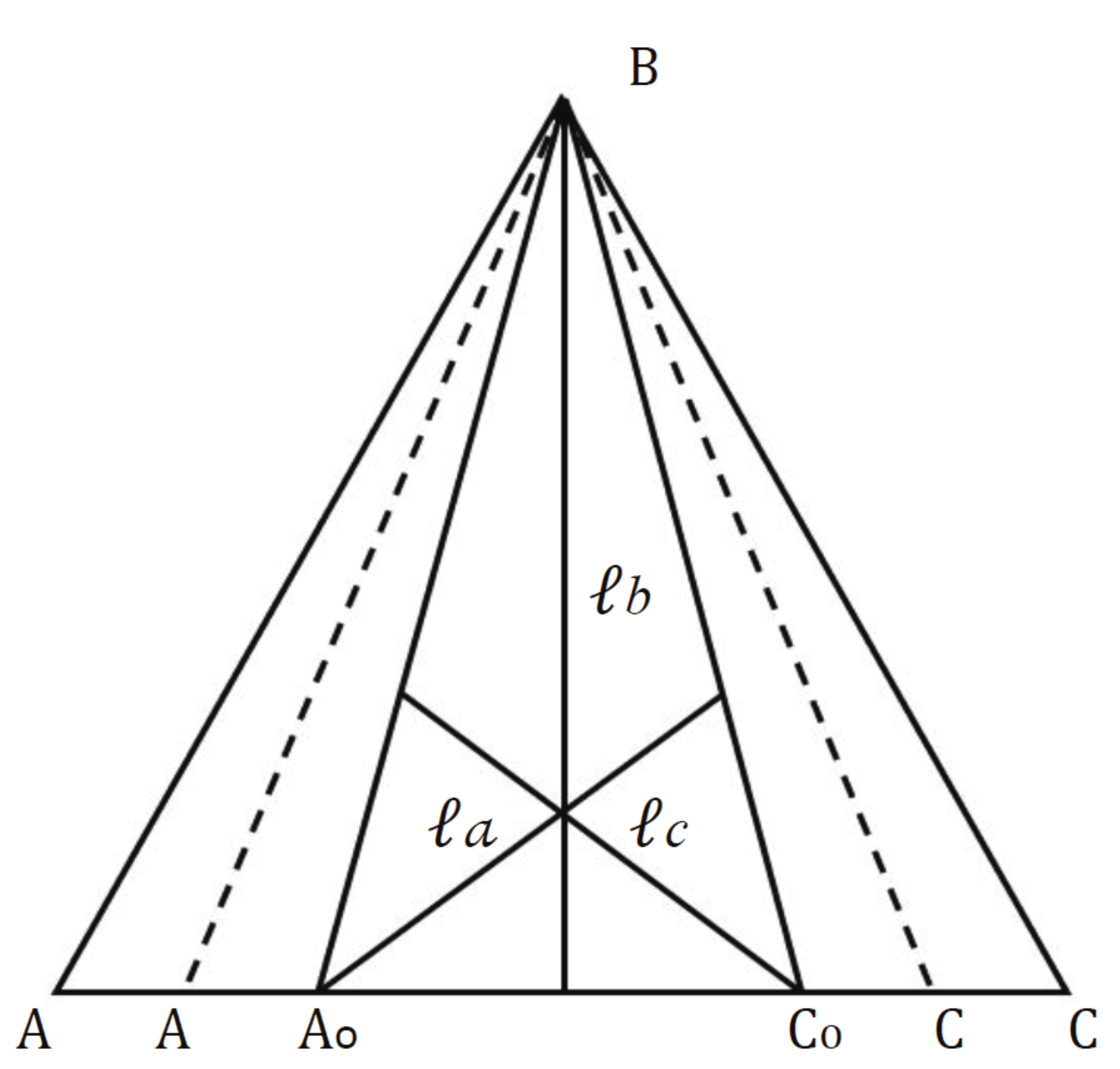}
\caption{Transformation of equilateral triangle $ABC$ with bisectors $l_a=l_b=l_c=l_1$ to isosceles triangle
$A_0BC_0$ with bisectors $l_a=l_c=l_2$ and $l_b=l_1$ by decreasing the angle $\angle B$.
 }
\end{center}\label{fig2}
\end{figure}

The above informal description of the idea of the proof actually introduces the
elements of the triangle $ABC$ as functions of the angle $\angle B$ under conditions
that lengths of one or two bisectors are kept constant. Therefore
it is convenient to distinguish the changing variables $l_a,\,l_b,\,l_c$
equal to the bisectors lengths of our varying triangle from their prescribed
fixed values. To this end, from now on we denote the fixed prescribed values as
$l_3<l_2<l_1$.

We shall start with formulation of an obvious Lemma 1 that describes the variation
of elements of an isosceles triangle in the first type of our transformations.

\smallskip

{\bf Lemma 1.} In an isosceles triangle $ABC$ with equal sides $AB=CB$ and fixed
length of the bisector $l_b=l_1$ (see Fig.~1) the angles
$\angle A = \angle B= (180^{\circ}-\angle B)/2$  are decreasing functions of $\angle B$
and the sides $AB=CB$ as well as bisectors $l_a=l_c$ are increasing functions of
$\angle B$. When $\angle B$ changes from $60^{\circ}$ to zero, the bisectors $l_a,\,l_c$
change from $l_b=l_1$ to zero. $\square$

\smallskip

In the second type of our transformations the length of two bisectors $l_b$ and $l_a$ are
kept constant and we are interested in dependence of the angles $\angle A$, $\angle C$
and the length of the bisector $l_c$ on the angle $\angle B$.

\smallskip

{\bf Lemma 2.} Let in a triangle $ABC$ the bisector $l_a$ (of the angle $\angle A$) be equal to
$l_a=l_2$, the bisector $l_b$ (of angle $\angle B$) be equal to $l_b=l_1$, and the angles $\angle A$, $\angle B$, $\angle C$ satisfy the condition $\angle B<\angle A< \angle C$ (see Fig.~2).
Then the following inequalities hold
\begin{equation}\label{eq1}
    \frac12l_2< b< a< c < l_1+l_2.
\end{equation}

\begin{figure}[th]
\begin{center}
\includegraphics[width=6cm,height=7cm,clip]{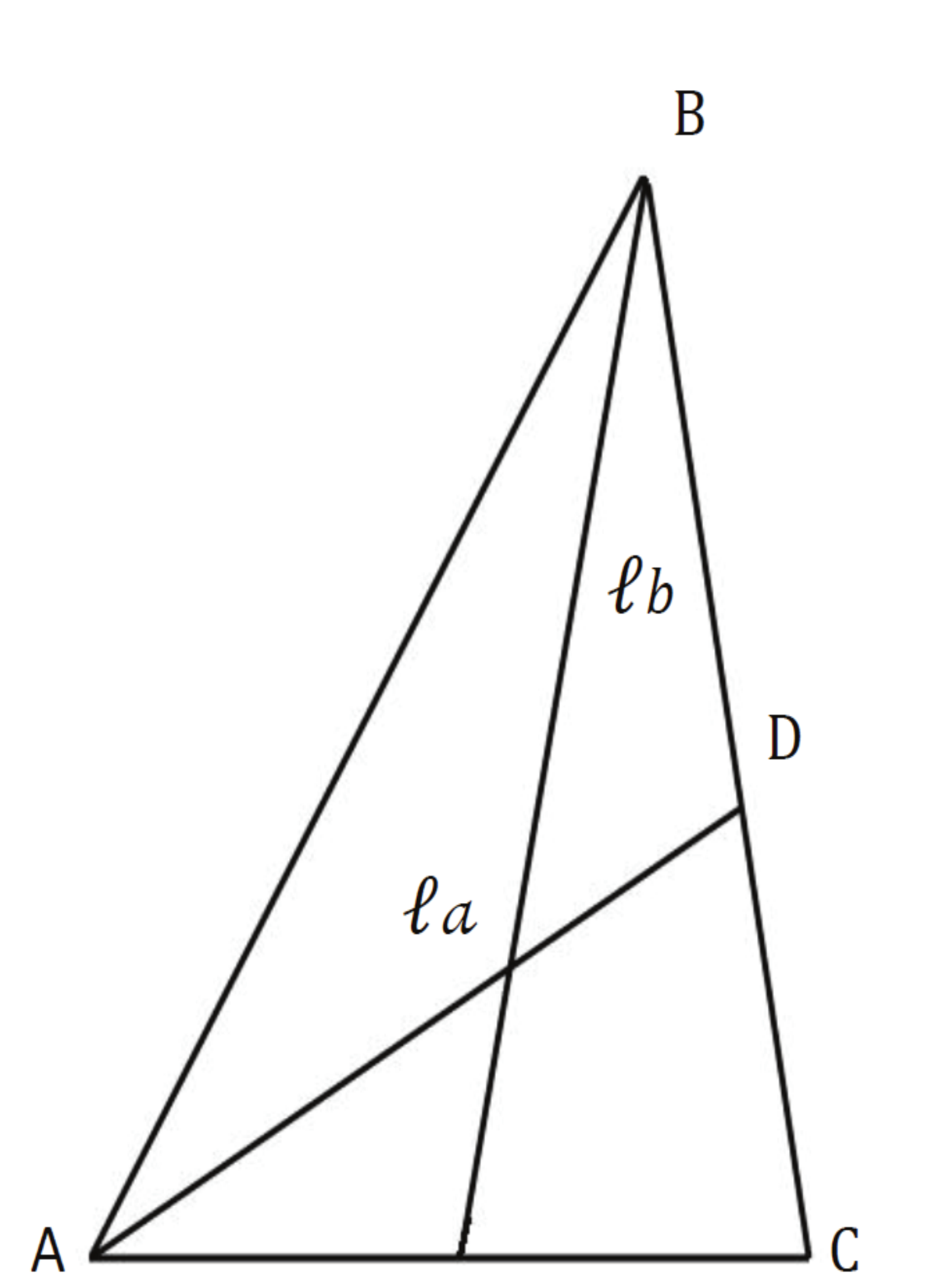}
\caption{Sketch of triangle $\angle B<\angle A< \angle C$ and $l_a=l_2< l_b=l_1$.
 }
\end{center}\label{fig1}
\end{figure}

{\bf Proof.}
Since $\angle ADC = \angle BAD+\angle ABD$, we have $\angle ADC>\angle BAD=\angle CAD$ and,
consequently, $CD<AC$. Therefore, we have $AD<AC+CD<2AC$, hence
$AC = b>\frac12AD=\frac12 l_2$. At last, $AB=c$ is less than the sum of the two other sides of
a triangle formed by $AB$ and pieces of the bisectors $l_b=l_1$
and $l_a=l_2$. $\square$

\smallskip

Now we consider a continuous set of triangles $ABC$ parameterized by their smallest
angle $\angle B$,
so that $\angle B<\angle A< \angle C$, and we assume that in all these triangles
the bisectors $l_a$ and $l_b$ have the same values. Thus, the angles $\angle A$ and $\angle C$
become some functions of the angle $\angle B$: if we change the angle $\angle B$,
then the angles $\angle A$ and $\angle C$ will also change whereas the bisectors
$l_a$ and $l_b$ remain constant. At first we shall prove the following simple property
of these functions.

\smallskip

{\bf Lemma 3.} If there exists such a transformation of a triangle $ABC$ that the conditions
$\angle B<\angle A< \angle C$ are satisfied, the bisectors $l_a$ and $l_b$ are kept constant,
and $\angle B\to0$, then $\angle A\to 0$ and $\angle C\to 180^\circ$.

\smallskip

{\bf Proof.} From identities
$\sin\angle A/a=\sin\angle B/b=\sin\angle C/c$ and inequalities
(\ref{eq1}) we find that in the limit $\angle B\to 0^\circ$ and, consequently,
$\sin\angle B\to 0$, we have $\sin\angle A\to0^\circ$,
$\sin\angle C\to 0^\circ$. These limiting values of $\angle A$ and $\angle B$
can be realized in two cases
\begin{equation}\nonumber
    \begin{split}
    & \angle A\to0^\circ,\qquad \angle C\to 180^\circ;\\
    & \angle A\to180^\circ,\qquad \angle C\to 0^\circ.
    \end{split}
\end{equation}
However in the second case we must have $l_a\to0$ and that contradicts
to the conditions of the Lemma,
so this case must be excluded.
We note also that inequalities (\ref{eq1}) exclude such ``exotic'' situations
as $\angle B\to0^\circ$, $b\to0$
or $\angle B\to 0^\circ$, $a\to\infty$, $c\to\infty$ that admit any
limiting values of angles $A$ and $B$ provided
$\angle A+\angle C\to 180^\circ$ and $\angle A<\angle C$.
This completes the proof of Lemma 3.
$\square$

\smallskip

\begin{figure}[t]
\begin{center}
\includegraphics[width=6cm,height=7cm,clip]{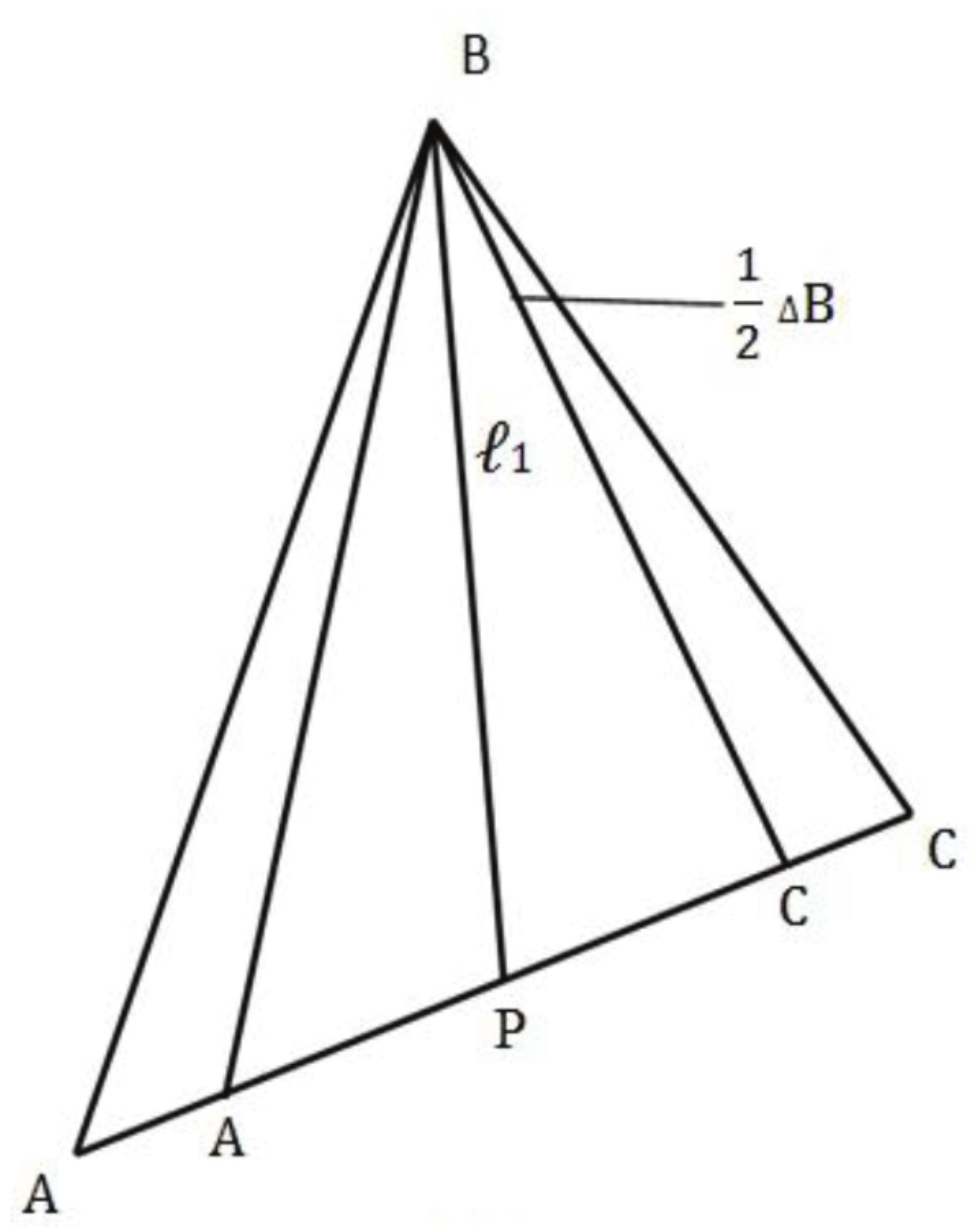}
\includegraphics[width=5.5cm,height=7cm,clip]{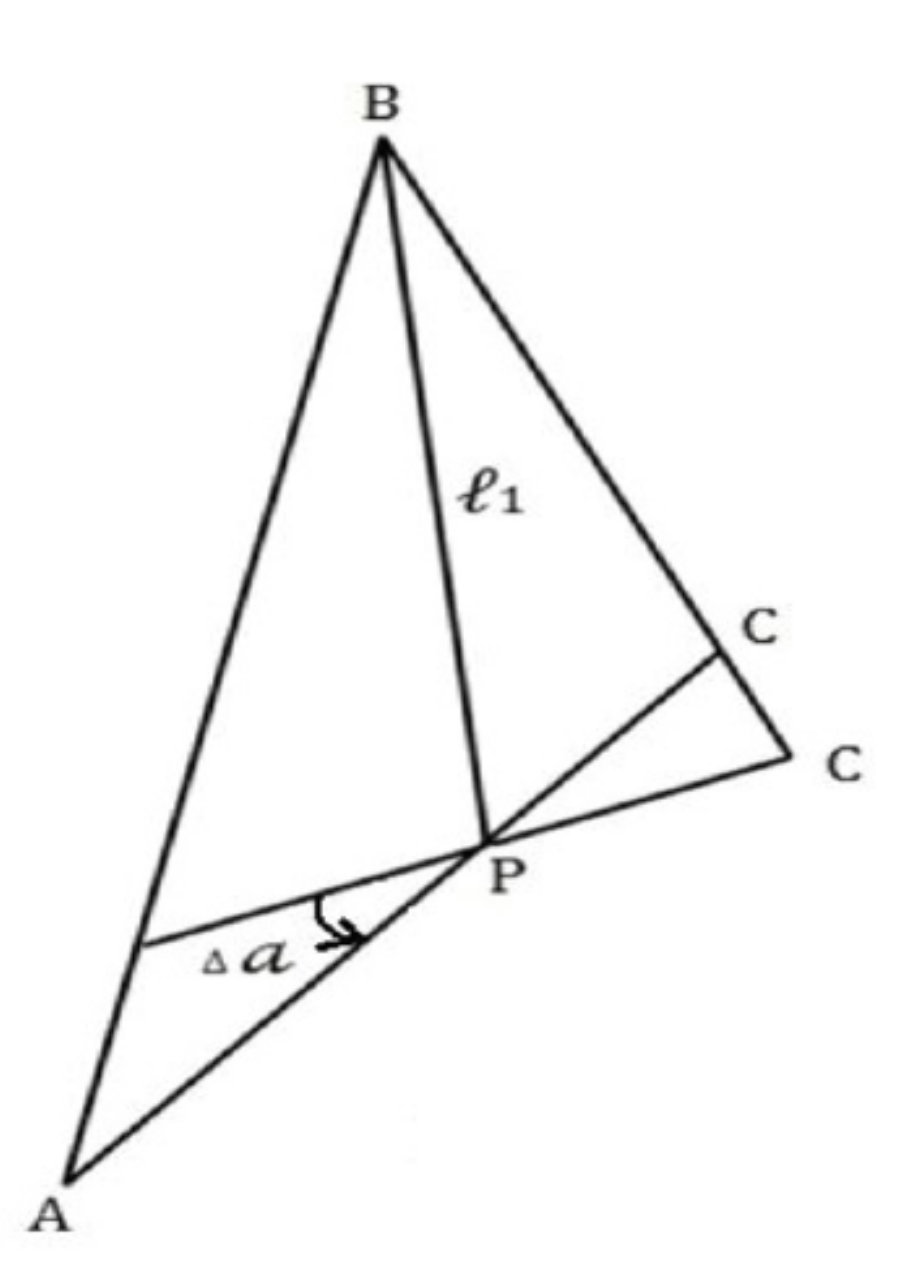}
\caption{Transformation of the isosceles triangle with bisectors
$l_a=l_c=l_2$ and $l_b=l_1$ to the final triangle
with bisectors with $l_b = l_1$, $l_a=l_2$ and $l_c=l_3$.
 }
\end{center}\label{fig3}
\end{figure}

Now we concretize the transformation of the second type by combining the decrease
of the angle $\angle B$ with rotation of the side $AC$ around the point $P$
at which the bisector $l_b$ and the side $AC$ meet. Let the angle $\angle B$ be
decreased by the value $\Delta B$ (see Fig.3, left) and the side $AC$ be rotated
counterclockwise around point $P$ by the angle $\Delta\alpha$ (Fig.~3, right).
Then after such a transformation
we have $\angle B<\angle C$ since $l_a<l_b$ and also we have $\angle A<\angle C$
since $\Delta A=\tfrac12\Delta B-\Delta \alpha$ and
$\Delta C=\tfrac12\Delta B+\Delta\alpha$. Thus, we have proved the following

\smallskip

{\bf Lemma 4.} In the defined above transformation we always have
$\angle B<\angle A<\angle C$. $\square$

\smallskip
This means that the conditions of Lemma 3 are always fulfilled for our
concrete realization of the transformation.

\smallskip

Now we shall study the behavior of the bisector $l_c$ under action of
our transformation.

\smallskip

{\bf Lemma 5.} When angle $\angle B$ decreases from its value $B_0$,
corresponding to the isosceles triangle with $l_a=l_c=l_2$, to zero, and
the bisectors $l_b=l_1$, $l_a=l_2$ are kept constant, then the length
of the bisector $l_c$ is a continuous function of $\angle B$ taking
all values from the interval $[0,l_2]$.

\smallskip

{\bf Proof.} From inspection of Fig.~3 (left) it is clear that decrease of
the angle $\angle B$ leads to the decrease of the sides $AB$ and $AC$
and to the increase of the angle $\angle A$. Consequently, from the known equation
\begin{equation}\label{eq2}
    l_a=\frac{2AB\cdot AC}{AB+AC}\cos\frac{\angle A}2=
    \frac2{\frac1{AB}+\frac1{AC}}\cos\frac{\angle A}2
\end{equation}
we find that the length $l_a$ also decreases with decrease of the angle $\angle B$,
when the angles $\angle A$ and $\angle C$ are kept constant.
To prevent such a decrease of $l_a$ with decrease of $\angle B$ by the value $\Delta B$,
we combine it with rotation of the side $AB$ around the point $P$ by the angle
$\Delta \alpha$, where $\alpha=\angle APB$.
Considered separately, such a rotation increases the sides $AB$ and $AC$ monotonously
to infinitely large values, decreasing at the same time the angle
$\angle A$ to infinitesimally small values. Then, according to Eq.~(\ref{eq2}),
the length $l_a$ can be increased by such a rotation as we
wish and, consequently, we can compensate the decrease of $l_a$ after diminishing of
$\angle B$ by the rotation of $AC$ and make the bisector of the angle $\angle A$
equal to $l_a=l_2$. This construction shows that each value of the angle $\angle B$
corresponds to a definite value of the angle $\angle A$ at which $l_a=l_2$ and
the function $\angle A=\angle A(\angle B)$ is continuous. The function $l_c(\angle B)$
is also continuous and in the described above transformation the conditions
$\angle B<\angle A< \angle C$, $l_a=l_2=\mathrm{const}$, $l_b=l_1=\mathrm{const}$
of Lemma 3 are fulfilled and when $\angle B$ tends to zero we have
$\angle A\to 0$ and $\angle C\to 180^\circ$, consequently $l_c\to0$ as $\angle B$
tends to zero.
Thus, we arrive at the conclusion that with decrease of $\angle B$ from $B_0$ to zero
the bisector $l_c$ as a continuous function of $\angle B$
takes all values from the interval $[0,l_2]$.

\smallskip

{\bf Remark.} As we shall see later, the function is not only continuous but also
single-valued and, hence, monotonous.

\smallskip

Now we can prove the main theorem about existence of a triangle with prescribed
lengths of its bisectors.

\smallskip

{\it Proof of the theorem.} ---
Let three segments be given whose lengths $l_1$, $l_2$, $l_3$ satisfy the conditions
$l_3<l_2<l_1$. Let us consider
an equilateral triangle $ABC$ with bisectors equal to $l_a=l_b=l_c=l_1$.
We fix the position and the length of the bisector
$l_b$ and decrease the angle $\angle B$ in a way shown in Fig.~1 down to the moment when we get an isosceles triangle
$A_0BC_0$ with bisectors $l_a=l_c=l_2$ and $l_b=l_1$. Existence of such a triangle
follows from Lemma 1: since the bisector lengths $l_a$ and $l_c$ ($l_a=l_c$) are
continuous functions of the angle $\angle B$ and they decrease
with $\angle B$ to zero value in the limit $\angle B\to0$, then they take
all values in the interval $[0,l_1]$
including the value $l_a=l_c=l_2$ that belongs to this interval due to inequality $0<l_2<l_1$.
The corresponding value of the angle $\angle B$ will be denoted as $\angle B_0$.

Then we continue transformation of the triangle $ABC$ but now keeping the values
of the bisectors $l_a$ and $l_b$ constant
and decreasing the length $l_c$ of the third bisector by further decrease of
the angle $\angle B$ with simultaneously rotation of
the side $AC$ counterclockwise around the point $P$ as is shown in Fig.~3.
As was shown in Lemma 5,
during such an evolution the function $l_c(\angle B)$  behaves as a
continuous function taking all values from the interval $[0,l_2]$
including the value $l_c=l_3$ thanks to inequality $0<l_3<l_2$. Hence
by means of our transformation we arrive at such a triangle whose bisectors have the
prescribed lengths.
Thus the main theorem is proved. $\square$

\smallskip

Now we shall make several additional remarks.

\smallskip

{\bf Remark.} In \cite{za-2003} an elementary proof is given that if a triangle
with prescribed lengths does exist then it is unique up to isometry. Combining this
theorem with our proof of existence of such a triangle, we conclude that in our
transformation a single value of each parameter $\angle A$, $\angle C$ and $l_c$
corresponds to a given value of $\angle B$; otherwise there would exist several such triangles.
Consequently, these continuous functions of the angle $\angle B$ must be
monotonous: with decrease of $\angle B$ the angle $\angle A$ and the bisector $l_c$
monotonously decrease whereas the angle $\angle C$ monotonously increases.

\smallskip

{\bf Remark.} The above proof shows that the angles of the final triangle satisfy
the inequalities $0^\circ<\angle A<\angle A_0$,
$0^\circ<\angle B<\angle B_0$, $\angle C_0<\angle C<180^\circ$,
where $\angle A_0=\angle C_0=(180^\circ-\angle B_0)/2$.
These inequalities can be made stronger if one introduces an angle $\beta$
such that $\angle C=180^\circ-2\beta$.
Let $\beta=\arcsin(l_3/(l_1+l_3))$, then we obtain the following estimate for $l_c$:
\begin{equation}\label{eq3}
    l_c=\frac{2ab}{a+b}\cos\frac{\angle C}2=\frac{2ab}{a+b}\cdot\frac{l_3}{l_1+l_3}<\frac{2a^2}{2a}\cdot\frac{l_3}{l_1+l_3}
    =l_3\cdot\frac{a}{l_1+l_3}<l_3.
\end{equation}
Since the angle $\angle C$ increases monotonously with decrease of $\angle B$
($\Delta C=\frac12\Delta B+\Delta\alpha$)
and the length $l_c$ monotonously decreases, then we find from
Eq.~(\ref{eq3}) that $\angle C<180^\circ-2\beta$.

\begin{figure}[th]
\begin{center}
\includegraphics[width=6cm,height=7cm,clip]{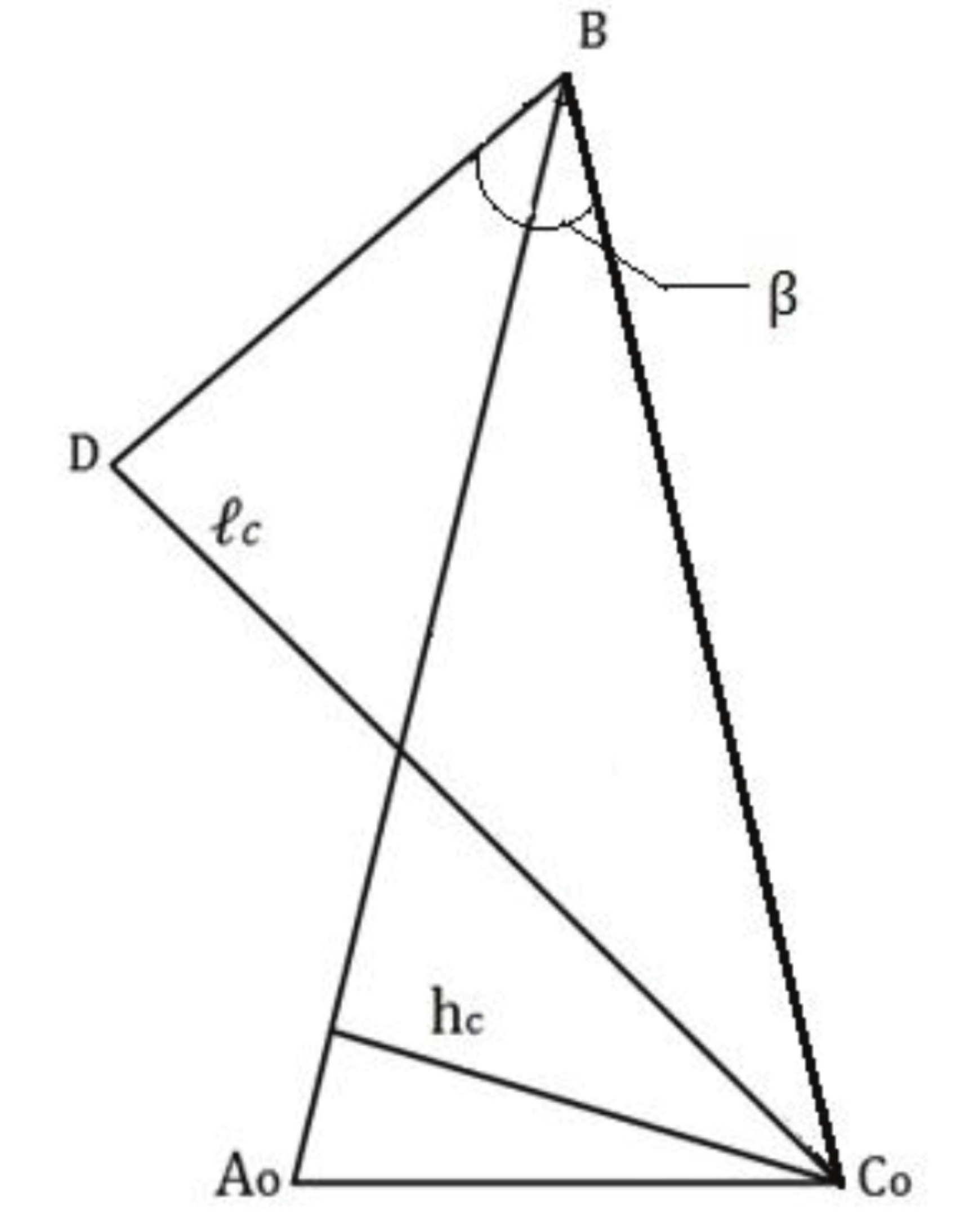}
\caption{Estimates of the angles $\angle B$ and $\angle A$.
 }
\end{center}\label{fig5}
\end{figure}

For estimation of the angle $\angle A$ we prolong the bisector line $l_c$
up to the point $D$ with $BD\bot CD$;
see Fig.~4. Then $\angle DBC=90^\circ-\frac12\angle C=\beta$ and we
get $\angle B<\beta$. Indeed, if we assume
that during our transformation at some moment we reach $\angle B>\beta$,
then by continuity argument there
exists a triangle $A'BC'$ such that $\angle B=\beta$. Consequently,
$\angle B+\frac12\angle C=\angle B+(90^\circ-\beta)=90^\circ$,
that is $AB\bot l_c$ and the triangle $A'BC'$ must be isosceles,
$AC=BC$, what is impossible because of $l_a\neq l_b$.
Hence, taking into account that $\angle A+\angle B+\angle C=180^\circ$,
we arrive at inequality $\angle A>\beta$.
Thus, with account of monotonicity of our functions, we conclude that
the angles of the triangle $ABC$ satisfy the inequalities
$\beta<\angle A<\angle A_0$,
$0^\circ<\angle B<\angle B_0$, $\angle C_0<\angle C<180^\circ-2\beta$.

\smallskip

{\bf Remark.} At last we notice that the suggested here method gives
actually the algorithm for numerical calculation of the elements of
the triangle with prescribed values of its bisectors.

\end{document}